\documentclass{amsart}
\usepackage{epsfig}
\newcommand{\lcr}{\raisebox{-5pt}{\mbox{}\hspace{1pt}
                  \epsfig{file=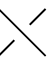}\hspace{1pt}\mbox{}}}
\newcommand{\ift}{\raisebox{-5pt}{\mbox{}\hspace{1pt}
                  \epsfig{file=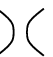}\hspace{1pt}\mbox{}}}
\newcommand{\zer}{\raisebox{-5pt}{\mbox{}\hspace{1pt}
                  \epsfig{file=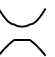}\hspace{1pt}\mbox{}}}

\newtheorem{theorem}{Theorem}[section]
\newtheorem{lemma}[theorem]{Lemma}
\newtheorem{prop}[theorem]{Proposition}

\theoremstyle{remark}
\newtheorem{remark}[theorem]{Remark}

\numberwithin{equation}{section}

\title[Noncommutative A-polynomial of figure 8]{The computation of
the non-commutative
generalization of the A-polynomial for  the figure-eight knot}
\author{R{\u{a}}zvan Gelca}
\address{Department of Mathematics and Statistics, 
Texas Tech University, Lubbock, TX 79409 and Institute of Mathematics
of the Romanian Academy, Bucharest, Romania}
\email{rgelca@math.ttu.edu}
\author{Jeremy Sain}
\address{Department of Mathematics, University of California at Berkeley,
Berkeley, CA 94720}
\email{jsain@Math.Berkeley.EDU}

\date{1 December 2002}


\newtheorem{cor}[theorem]{Corollary}

\begin{document}
\maketitle

\begin{abstract}
The paper computes the noncommutative A-ideal of
the figure-eight knot,  a noncommutative
generalization of the A-polynomial. It is  shown that if
a knot has the same A-ideal as the figure-eight knot,
then all colored Kauffman brackets are the same as those of the
figure eight knot.
\end{abstract}

\section{Introduction}

The noncommutative A-ideal is a knot invariant that generalizes
the A-polynomial of Cooper, Culler, Gillet, Long, and Shalen.
It was introduced in the framework of quantum topology.
Recall that the A-polynomial is computed
 from the  pull
 back  to the ring of polynomials,  through a sequence of rational
maps, of  the ideal that determines 
the character variety of the knot complement. The noncommutative A-ideal 
 is a deformation of this ideal with respect to
a parameter, and is  a finitely generated
 ideal   in the quantum plane.

The purpose of this paper is to find the noncommutative A-ideal
of the figure-eight knot. The computation is based on 
skein operations in the knot complement  and on the noncommutative
trigonometry of the noncommutative torus. While the formulas we
obtain are complicated and not at all illuminating, 
there are three observations that can be extracted from
this work.

 First, the fact that the noncommutative A-ideal of the 
figure-eight knot can be computed. Up to now only the A-ideals
of the unknot, trefoil, and to some extent $(2,2p+1)$-torus knots
have been determined. But the complements of all these knots are 
Seifert fibered spaces, which is a true advantage when computing 
 quantum invariants.  The present paper is the first example
of a computation of the noncommutative A-ideal for a hyperbolic
knot. 

Second, our methods can be adapted to give a pictorial way
of computing the (classical) A-polynomial of the figure-eight knot.
The work is  done in the ring of functions on the
$SL(2,{\mathbb C})$-character variety of the knot complement.
It  will be seen below that  an important role is played
by the characters of irreducible representations of $SL(2,{\mathbb C})$.

Finally, as a corollary of our work, we discovered that the
noncommutative A-ideal of the figure-eight knot determines
all  colored Kauffman brackets of this knot. More precisely, 
any knot that has the same A-ideal as the figure-eight knot
has the same colored Kauffman brackets, in the zero framing of the knot.
Recall that the colored Kauffman brackets are obtained by
evaluating  the classical Kauffman bracket on the colorings of
 the knot by Jones-Wenzl idempotents. 

This property is a consequence of the so called ``orthogonality relation''
between the A-ideal and the Jones polynomial. 
It has been proved  before that the noncommutative A-ideals of the
unknot and the $(2,2p+1)$-torus knots determine their colored Kauffman
brackets. We now  suggest the hypothesis  that this property  
is true for  all knots.

The paper is organized as follows. In Section 2 we recall the definitions
of the A-polynomial and of the noncommutative A-ideal, together with some
basic facts about the structure of the Kauffman bracket skein algebra
of the torus. In Section 3 we explain briefly how the Kauffman bracket
skein module of the complement of the figure-eight knot is computed,
which is the content of an unpublished  theorem of Bullock and Lofaro.
In Section 4 we perform the skein computations that carry all
specific information about the action of the Kauffman bracket skein
algebra of the torus on the skein module of the knot complement. 
Section 5 computes the generators of the peripheral ideal
of the knot, which is  the left ideal  sent to zero by the map
from the skein algebra of the boundary torus to the skein
module of the knot complement. This ideal determines the noncommutative
A-ideal, exhibited in the next section. Section 7 contains what
we consider the most important result of the paper, a theorem showing 
that the A-ideal of the figure-eight knot
 determines all its colored Kauffman brackets.

The authors would like to thank Fumikazu Nagasato for 
discovering and correcting an error that appeared in Proposition 5.4 below.

\section{Some facts about the A-polynomial and the noncommutative A-ideal}

The computation of  the A-polynomial of a knot  is based on the  
following sequence of rational maps
\begin{eqnarray*}
{\mathbb C}\times {\mathbb C}\hookleftarrow {\mathbb C}^*\times
{\mathbb C}^*\rightarrow X({\mathbb T}^2)\leftarrow X(M),
\end{eqnarray*}
where $X(M)$ and $X({\mathbb T}^2)$ are the character varieties
of $SL(2,{\mathbb C})$-representations of the fundamental 
group of the knot complement, respectively the torus. The
first map is the inclusion, the second is the 2-1 covering map
of the character variety, and the third  is obtain by 
restricting representations to loops on the boundary of the knot
complement. If we send
$X(M)$ to ${\mathbb C}\times {\mathbb C}$ through this
sequence,  we obtain a closed  algebraic
set whose one dimensional part is determined by the A-polynomial (times
$l-1$) (see \cite{CCGLS}).

Considering the rings of regular functions on these varieties,
we obtain the dual sequence of maps
\begin{eqnarray*}
{\mathbb C}[l,m]\hookrightarrow {\mathbb C}[l,l^{-1},m,m^{-1}]
\hookleftarrow {\mathbb C}[l,l^{-1},m,m^{-1}]^{inv} \rightarrow {\mathcal F}(
X(M))
\end{eqnarray*}
where $ {\mathbb C}[l,l^{-1},m,m^{-1}]^{inv}={\mathcal F}(X({\mathbb T}^2))
$ is the ring of functions on the character variety
of the torus, which  consists of  the Laurent polynomials
that are 
invariant under the simultaneous transformations $l\rightarrow l^{-1}$, and
$m\rightarrow m^{-1}$; and ${\mathcal F}(X(M))$ is the ring of
 regular functions
on the character variety of the knot complement. 
Let $I(K)\in  {\mathbb C}[l,l^{-1},m,m^{-1}]^{inv}$ be the 
kernel of the third map of this sequence. Extending $I(K)$
to the ring of Laurent polynomials and then contracting it to
${\mathbb C}[l,m]$ yields   a finitely generated
ideal of polynomials. The radical of the one dimensional part of this
ideal is generated by the A-polynomial (times $l-1$).

Now recall that the Kauffman bracket skein module of a 3-manifold $M$,
denoted by $K_t(M)$,  is
the quotient  of the free ${\mathbb C}[t,t^{-1}]$-module  with basis
the isotopy classes of links in $M$ by the submodule spanned by
relations of the form  $\displaystyle{\lcr-t\zer-t^{-1}\ift}$
and 
$\bigcirc+t^2+t^{-2}$, where the links in each expression are
identical except in a ball in which they look like depicted.
Since in the computations below we will divide by polynomials,
we slightly change the point of view, working over the field
${\mathbb C}(t)$ instead of the ring of Laurent polynomials
${\mathbb C}[t,t^{-1}]$. Note that in \cite{frgelo} $t$ was 
a complex number instead of a variable; all considerations
below work in this setting, too,  when  $t$ is not a $48$th root of unity.  
So our modules are in fact vector spaces, but for historical considerations
we will keep the name. 

At $t={-1}$, $K_{-1}(M)$ endowed with
the multiplication induced by  the disjoint union of links, is a ring. 
As explained in
\cite{bullock2} and \cite{przsik},  when 
reduced by $\sqrt{0}$, it  is isomorphic to 
${\mathcal F}(X(M))$, the ring of regular function on the 
$SL(2,{\mathbb C})$-character variety.  

For the cylinder over a surface $\Sigma$,
$K_t(\Sigma\times I)$ is an algebra. Here and below $I$ is 
the interval $[0,1]$. 
It was shown in \cite{FG} that the Kauffman
bracket skein algebra of the cylinder over the  torus is isomorphic to
the algebra generated by noncommutative cosines in the noncommutative
torus. The proof relies  on the product-to-sum formula
in $K_t({\mathbb T}^2\times I)$, which will play an important role
in this paper,  so we remind it to the reader. Let $T_n(x)$ be
the $n$th Chebyshev polynomial, $T_0(x)=2$, $T_1(x)=x$,
$T_{n+1}(x)=xT_n(x)-T_{n-1}(x)$, $n\in {\mathbb Z}$.
For two integers $p$ and $q$ let $n$ be their greatest common divisor,
$p'=p/n$, $q'=q/n$. The skein  $(p,q)_T$ is by definition 
equal to $T_n((p',q'))$
where the variable of the Chebyshev polynomial is the $p'/q'$ curve
on the torus, and powers consist of parallel copies. Then
$\{(p,q)_T, \; p,q\in {\mathbb Z},p\geq 0\}$ is
a basis of $K_t({\mathbb T}^2\times I)$. The product-to-sum formula
is 
\begin{eqnarray*}
(p,q)_T*(r,s)_T=t^{|^{pq}_{rs}|}
(p+r, q+s)_T+
t^{-|^{pq}_{rs}|}
(p-r,q-s)_T.
\end{eqnarray*}
The elements $t^{-pq}l^pm^q$ are the noncommutative exponentials of
the noncommutative torus, so the map 
$(p,q)_T\rightarrow t^{-pq}(l^pm^q+l^{-p}m^{-q})$ identifies 
 $\frac{1}{2}(p,q)_T$  with a noncommutative
cosine. 

 We  can now deform the  sequence of maps  defining
the A-polynomial to the following
\begin{eqnarray*}
{\mathbb C}_t[l,m]\hookrightarrow {\mathbb C}_t[l,l^{-1},m,m^{-1}]
\hookleftarrow K_t({\mathbb T}^2\times I) \rightarrow K_t(M).
\end{eqnarray*}
Here 
$ {\mathbb C}_t[l,m]$ and ${\mathbb C}_t[l,l^{-1},m,m^{-1}]$
are the rings of polynomials, respectively Laurent polynomials,
in the noncommuting variables $l$ and $m$, satisfying
$lm=t^2ml$. Let $I_t(K)$ be the kernel of the 
third map, which we  call the {\em peripheral ideal}. It is a left
ideal.  If we extend
it to Laurent polynomials, then contract it to $ {\mathbb C}_t[l,m]$,
we obtain a finitely generated left  ideal of polynomials, called
the noncommutative A-ideal of the knot. It was introduced in \cite{frgelo}
as a generalization of the A-polynomial. In this paper
we  compute it for the figure-eight knot.

At this point we would like to recall 
 that in quantum mechanics ${\mathbb C}_t[l,l^{-1},m,m^{-1}]$ is 
 the algebra  generated by the exponentials of the momentum and position
operators and their inverses. This is the ring of trigonometric 
polynomials in the noncommutative torus. 
The ring ${\mathbb C}_t[l,m]$ is called the quantum plane.

For people interested ultimately in the A-polynomial we
make the following remarks. The A-polynomial is usually determined by
looking first at the representation variety of the fundamental group,
then using elimination theory to obtain relations at the level of
character variety. In this sense the computation is ``geometric''
- it is done at the level of the underlying geometric space.
 
Now if in the computations of this paper
 we substitute $t$ by $-1$, then, with
appropriate modifications to account for the problem with 
 $48$th roots of unity (see below),
we obtain an ``algebraic'' computation of the A-polynomial. 
This means that  we study 
the ring of functions on the character variety, instead of the character 
variety itself. 
A nice feature of this approach is that it can be done pictorially.
A skein is  identified with  the function on
the character variety defined by  evaluating  characters
on that  skein \cite{CS}, \cite{bullock2},
\cite{przsik}. The sliding under a crossing corresponds to the
element in the ideal of the variety that comes from the
crossing, in the Wirtinger presentation of the knot. The Kauffman
bracket skein relation at $t=-1$ 
corresponds to the Hamilton-Cayley identity for $SL(2,{\mathbb C})$.
This method was  applied in 
\cite{naga} for computing the A-polynomial
of $(2,2p+1)$-torus knots. 

One should note that an important role in these computations
is played by the polynomials $S_n(x)$, where $S_0(x)=1$,
$S_1(x)=x$, $S_{n+1}(x)=xS_n(x)-S_{n-1}(x)$, $n\in {\mathbb Z}$.
For $n<0$, $S_n(x)=-S_{-n-2}(x)$, while for $n>0$, $S_n(x)$ is
the {\em character of the $n+1$-dimensional irreducible representation
of $SL(2,{\mathbb C})$}. These polynomials
 satisfy $S_n-S_{n-2}=T_n$, where
$T_n$ is the $n$th Chebyshev polynomial.

\section{The Kauffman bracket skein module of the complement
of the figure-eight knot}

It was shown in \cite{bullof} that the Kauffman bracket skein module of
the complement of the figure-eight knot is the ${\mathbb C}(t)$-module
with basis $x^n, x^ny,x^ny^2$, $n\geq 0$, or equivalently
$x^n$, $x^nz$, $x^nz^2$, $n\geq 0$, where $x$, $y$ and  $z$ are
depicted in Figure 3.1. Here and below, for two
curves $\gamma _1$  and $\gamma  _2$ when we write $\gamma _1^n\gamma _2^m$
we understand the disjoint union of $n$ parallel copies of the
first curve and $m$ parallel copies of the second, all copies
lying in  (very thin) tubular neighborhoods of the initial curves.

\begin{figure}[htbp]
\centering
\leavevmode
\epsfxsize=1.3in
\epsfysize=1.3in
\epsfbox{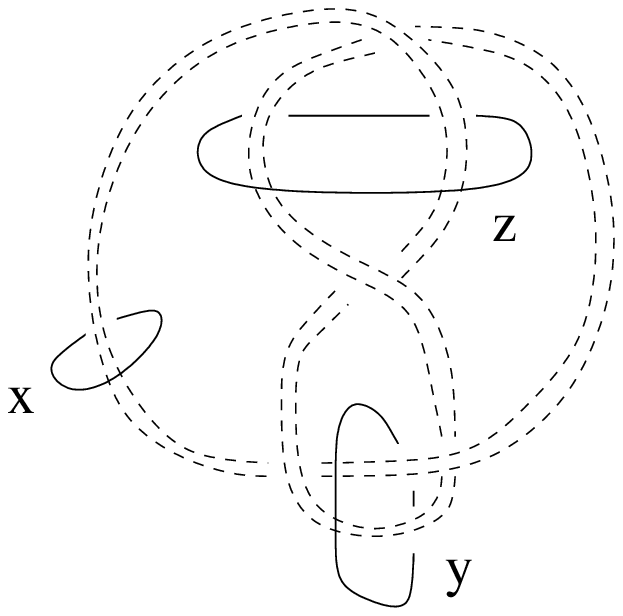}

Figure 3.1.    
\end{figure}

Since \cite{bullof} has not been published yet, we give here a 
brief overview of how  the Kauffman bracket skein module
of the complement of the figure-eight knot is computed. The complement of
this knot is obtained by attaching a 2-handle
to a genus two handlebody along the curve described
in Figure 3.2. In the picture one should understand
that the curve lies on the  sphere, on viewer's side.

\begin{figure}[htbp]
\centering
\leavevmode
\epsfxsize=1.8in
\epsfysize=1.3in
\epsfbox{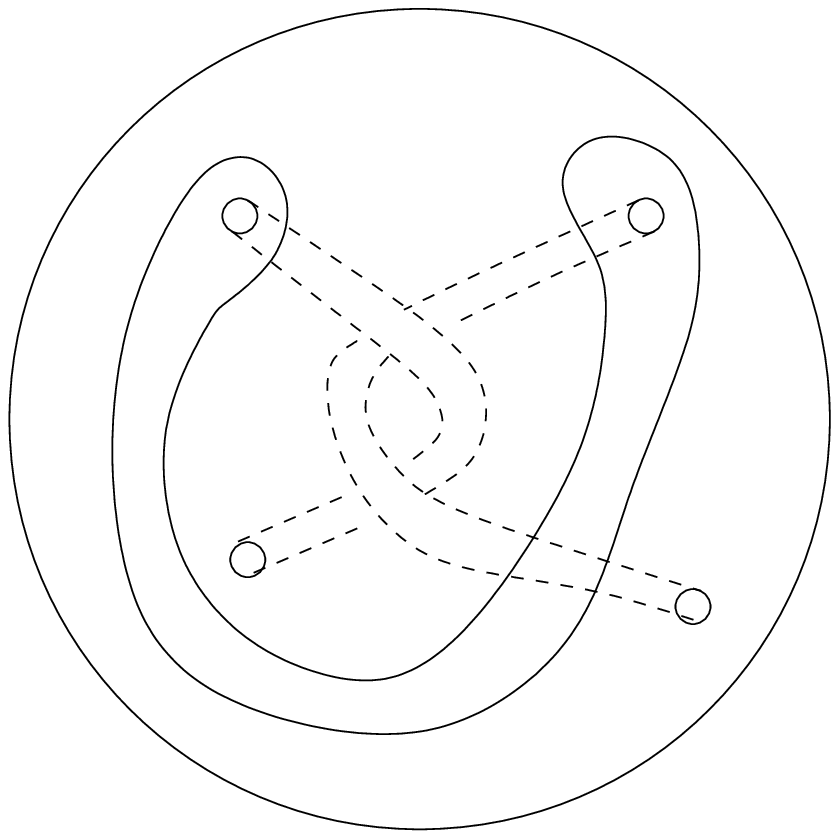}

Figure 3.2.    
\end{figure}

Represent the handlebody in the standard way, and push the curve 
inside slightly as in Figure 3.3. The two images represent the
side view and  the view
from above  of the handlebody. 

\begin{figure}[htbp]
\centering
\leavevmode
\epsfxsize=4in
\epsfysize=1.2in
\epsfbox{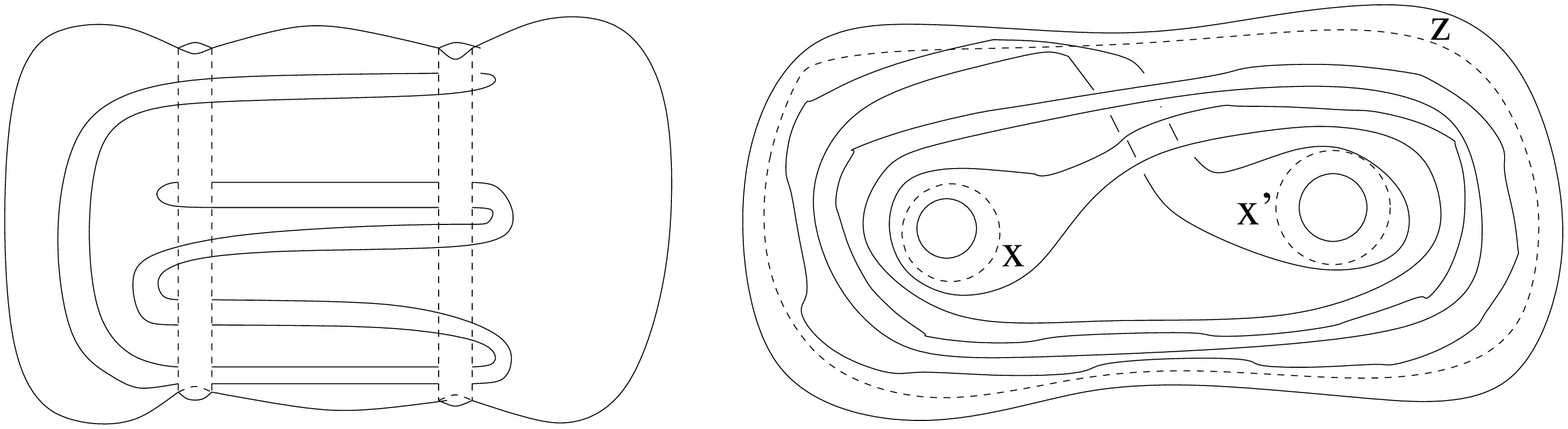}

Figure 3.3.    
\end{figure}

The Kauffman bracket skein module of the handlebody,
$K_t(H)$ is free with basis $x^nx'^mz^k, m,n,k\geq 0$, where 
$x,x',z$ are the dotted curves from  Figure 3.3 \cite{P}.
 The skein module of the knot
complement is obtained by factoring by the submodule   corresponding
to slidings of skeins through the 2-handle. This submodule 
is generated by  slidings  of  basis elements. 
Clearly after attaching the handle $x$ and $x'$ are isotopic.
Thus we conclude that the Kauffman bracket skein module of the 
knot complement is obtained by factoring the 
free module with basis $x^nz^m$ by the submodule $J$ generated
by $x^nz^m-x^nz^m\cdot sl(z)$, where $sl(z)$ is the slide of $z$ through
the 2-handle. 

Now sliding through the 2-handle amounts to taking
the banded sum with the attaching curve. Of course 
this can be done in infinitely many ways, but as explained
in \cite{bullock}, we only need to consider one such banded sum,
and we will always consider the one from 
Figure 3.4. 

 \begin{figure}[htbp]
\centering
\leavevmode
\epsfxsize=4in
\epsfysize=1.4in
\epsfbox{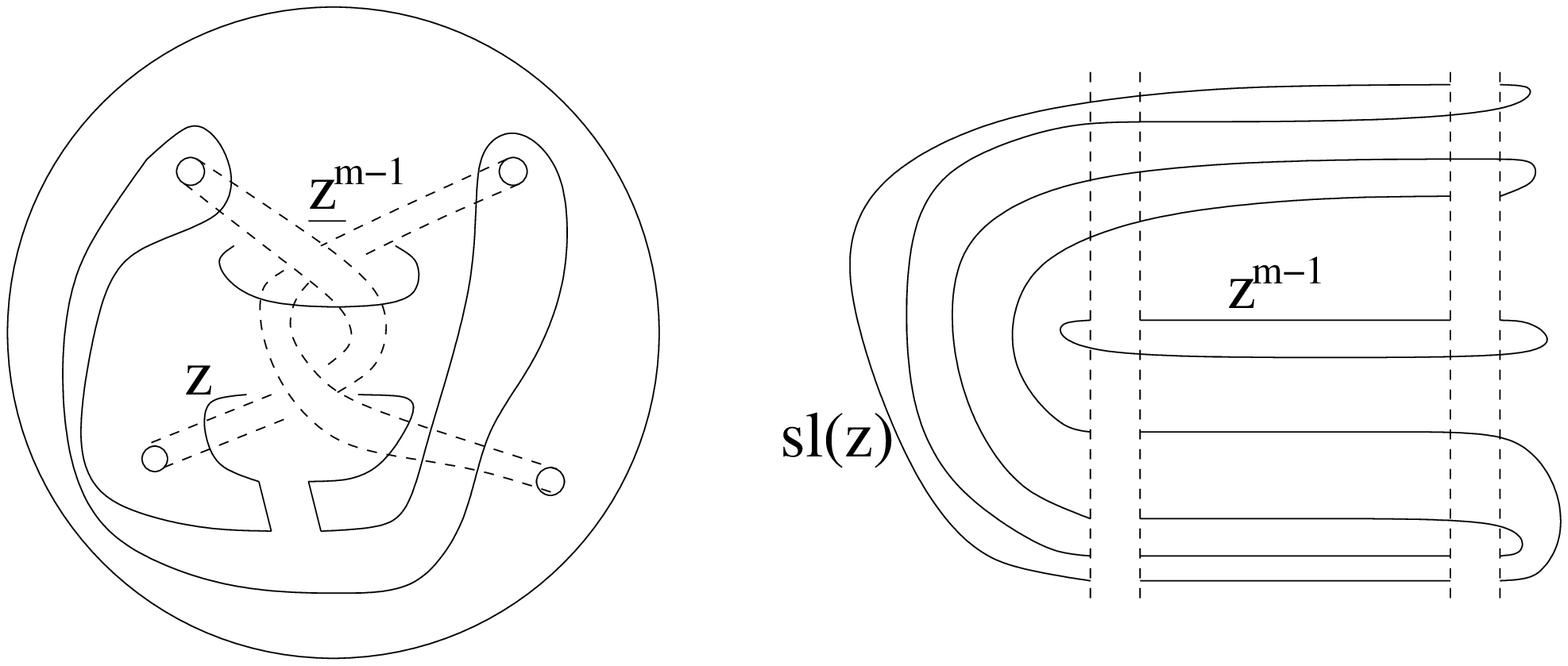}

Figure 3.4.    
\end{figure}

It is not hard to
 see that up to a power of $t$, the element $x^nz^m-x^nz^{m-1}
sl(z)$, $m\geq 1$  is a polynomial of $(m+3)$rd 
degree in $z$ with the coefficient of $z^{m+3}$ equal to $x^n$
(just count how many times the curves wind ``around'' the handlebody
and think what happens when you resolve the crossings).
When we slide the trivial  link (that is when we consider the 
attaching curve as a skein itself), we obtain a polynomial
of the form $z^4+$ terms of lower degree. 
Through an unpleasant but routine computation one finds
this polynomial, and  also the ones for $z-sl(z)$ and $z^2-zsl(z)$. Linear
combinations of them yield   polynomials of degrees $3$, $4$, and $5$ in $z$.
Hence $J$ is generated by polynomials of the form $x^nz^m+$
terms of lower degree in $z$, one for each $m\geq 3$. 
This implies that the Kauffman
bracket skein module of the knot complement is free with
basis $x^n$, $x^nz$, $x^nz^2$, $n\geq 0$. 

Lemma \ref{yz} below shows that we can replace this basis by
 $x^n$, $x^ny$, $x^nz$, $n\geq 0$, 
 where $x,y,z$ are the curves from Figure 3.1.
This is the basis that we will  work with.

We now make an {\em important observation} related to the amphichirality of
the figure-eight knot. If we take the mirror image,  the
knot maps to itself. Since $x$ is a meridian, it 
stays unchanged, but $y$  maps to $z$ and $z$ to $y$. 
Also, under this symmetry $t$ changes to $t^{-1}$ and the skein $(p,q)_T$
on the boundary changes to $(p,-q)_T$. 

\section{Skein computations in the knot complement}

\begin{lemma}\label{yz}
The following identities hold
\begin{alignat*}{1}
y^2&=-t^{-2}x^2y-t^{-6}z-2t^{-4}x^2+t^{-4}+1\\
z^2&=-t^2x^2z-t^6y-2t^4x^2+t^4+1.\\
\end{alignat*}
\end{lemma}

\begin{proof} 
%

To find $y^2$  start with the computation from Figure 4.1.
\begin{figure}[htbp]
\centering
\leavevmode
\epsfxsize=5in
\epsfysize=2.7in
\epsfbox{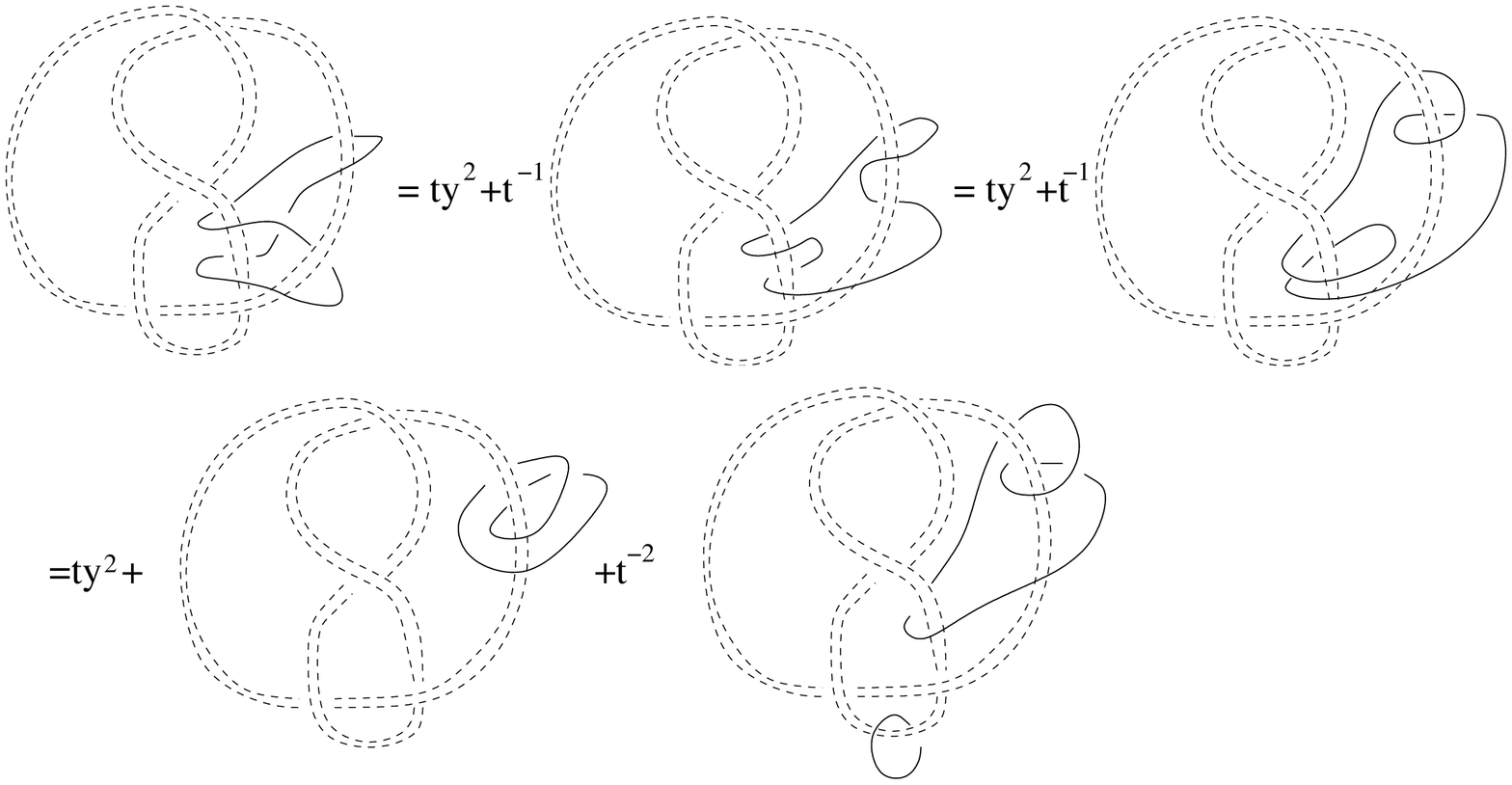}

Figure 4.1.    
\end{figure}
This expression is  equal to $ty^2+tx^2+t^{-1}(-t^2-t^{-2})+
t^{-1}x^2y+t^{-3}x^2$. 

On the other hand, performing an isotopy on the strand we started
with, we can then proceed as in Figure 4.2.
\begin{figure}[htbp]
\centering
\leavevmode
\epsfxsize=5in
\epsfysize=1.8in
\epsfbox{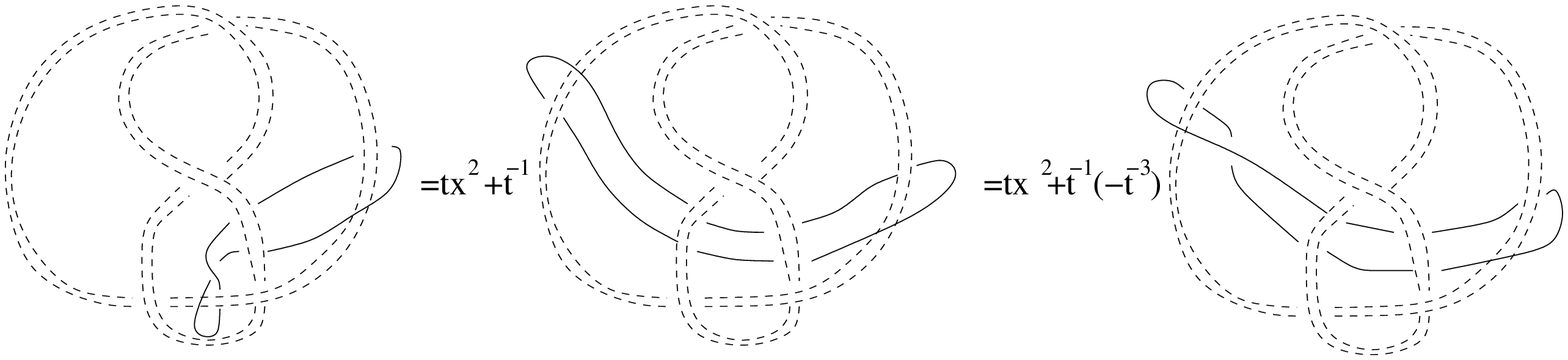}

Figure 4.2.    
\end{figure}
Resolving the crossing we obtain $-t^{-3}x^2+tx^2-t^{-5}z$.
The formula for $y^2$ now follows from the equality of these two expressions.
Take the mirror image to obtain the 
 formula for $z^2$.
\end{proof}

Let $M$ be the complement of a regular neighborhood of
the figure-eight knot complement. 
Below we will denote the image
of  $(p,q)_T\in K_t({\mathbb T}^2\times I)$ through the morphism
$K t({\mathbb T}^2\times I)\rightarrow K_t(M)$ also by  $(p,q)_T$.  

\begin{lemma}\label{(1,0)}
The images of $(1,0)_T$ and $(1,1)_T$ in the skein module of
the figure-eight knot complement are 
\begin{alignat*}{1}
(1,0)_T&=t^4x^2y+t^{-4}x^2z+(1-t^4)y+(1-t^{-4})z+2(t^2+t^{-2})x^2-t^2-t^{-2}\\
(1,1)_T&=t^5x^3y-2t^5xy+txz+t^{-3}xz+2t^3x^3-3t^3x+2t^{-1}x.
\end{alignat*}
\end{lemma}

\begin{proof}
We will only compute the image of $(1,0)_T$ in the skein module of the knot
complement, the image of $(1,1)_T$ is then found in an  analogous manner.
The computation is described in Figure 4.3, where we start by
resolving the two crossings on top. 
\begin{figure}[htbp]
\centering
\leavevmode
\epsfxsize=5in
\epsfysize=2.7in
\epsfbox{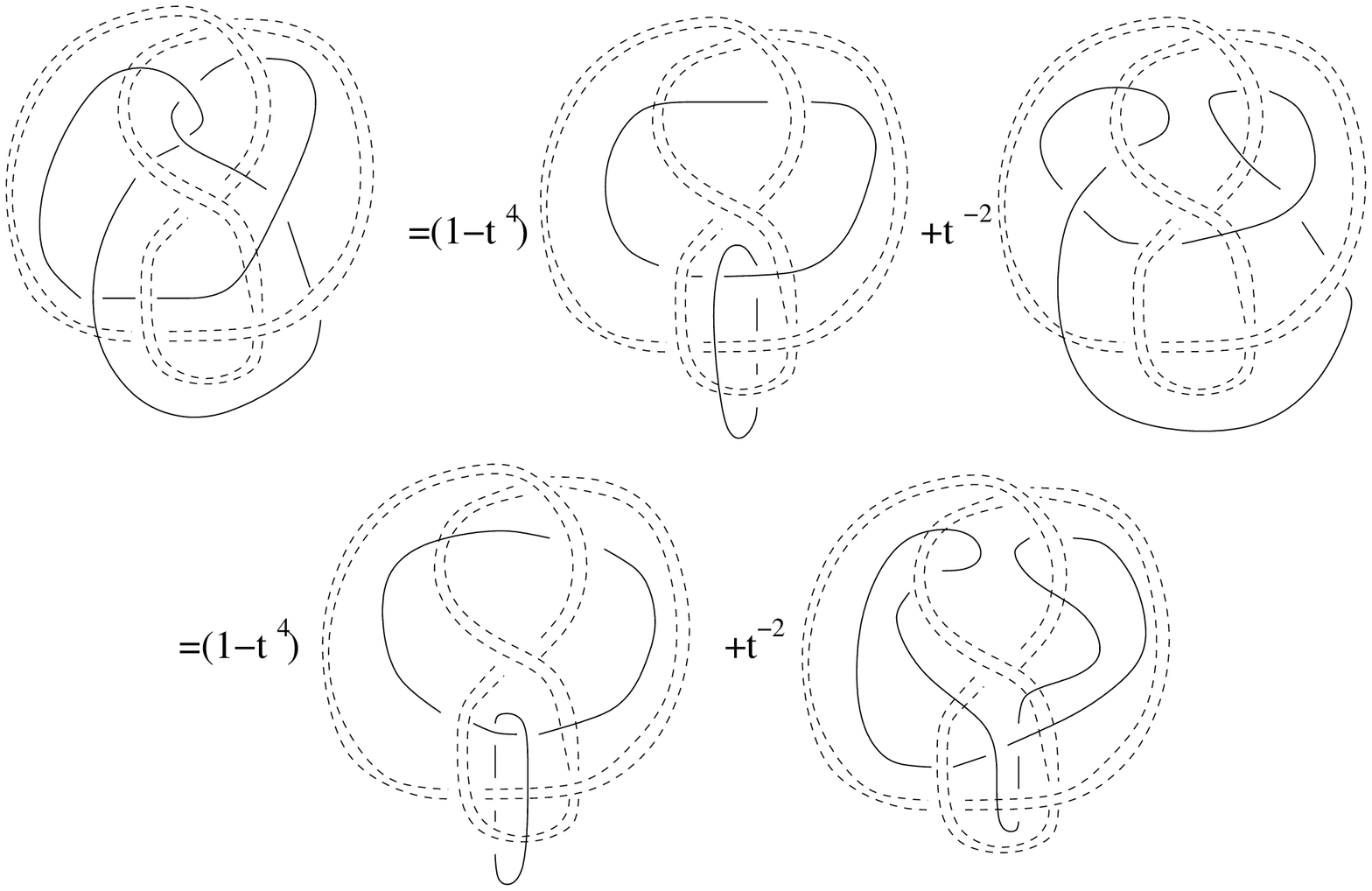}

Figure 4.3.    
\end{figure}

Call the two terms on the second line A and B. Resolving the two crossings
in A and performing the appropriate isotopies yields 
$-(1-t^4)(1-t^{-4})(t^2-t^{-2})+t^2(1-t^4)y^2$. We point out
the fact that the first term of this sum is obtained after replacing
a trivial knot by $-t^2-t^{-2}$. With B the story
is more complicated. We resolve the two crossings, then continue
as  in Figure 4.4.
\begin{figure}[htbp]
\centering
\leavevmode
\epsfxsize=5in
\epsfysize=3.5in
\epsfbox{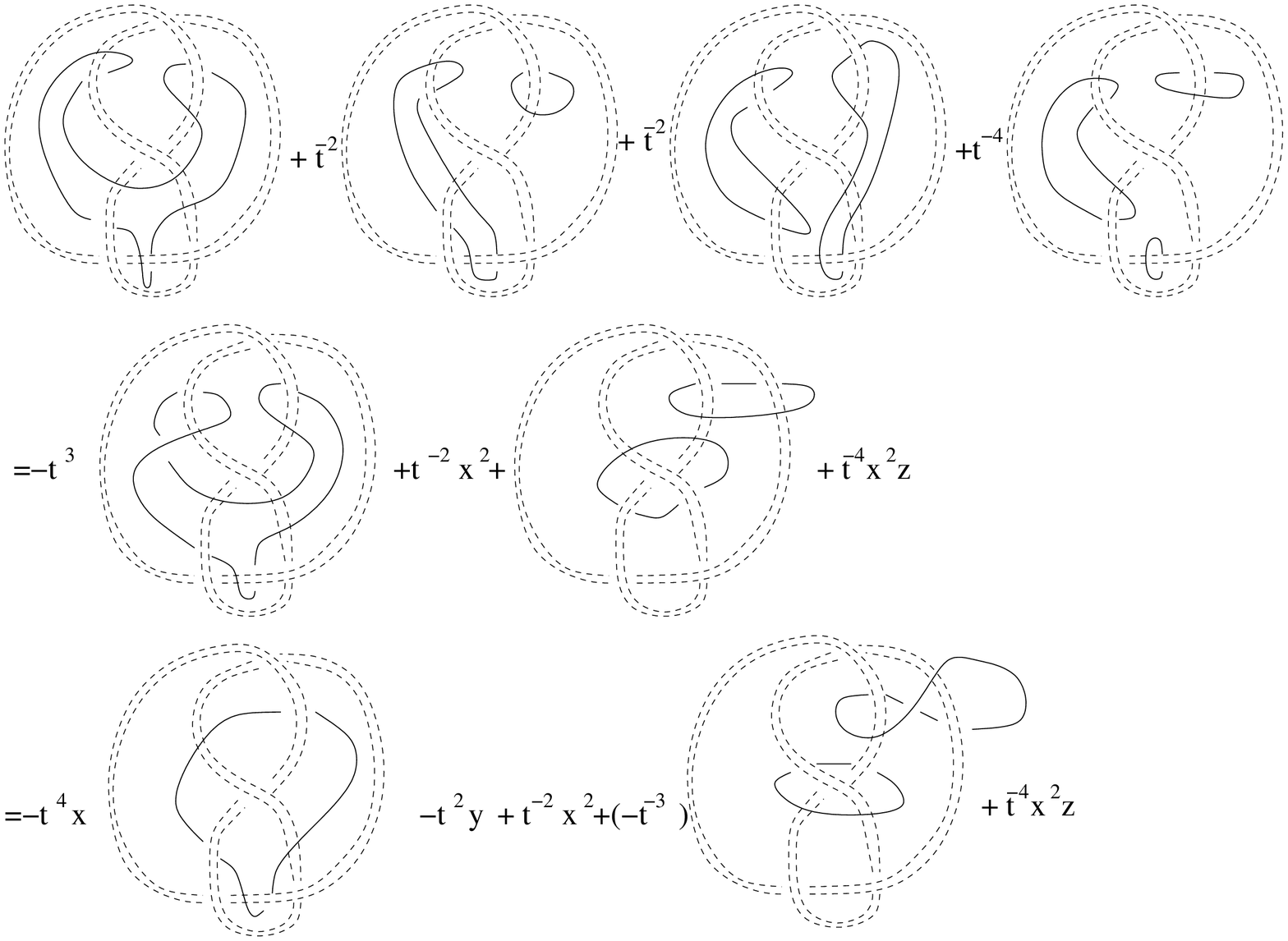}

Figure 4.4.    
\end{figure}

In the end, after resolving the crossing, the term 
involving the second diagram is equal to $-t^{-4}x^2z-t^{-6}z^2$.
In Figure 4.5 we 
compute the value of the remaining diagram (coefficient $-x$ 
excluded), which is equal to $-t^{-2}xy-t^{_4}x$. 
\begin{figure}[htbp]
\centering
\leavevmode
\epsfxsize=5in
\epsfysize=1.2in
\epsfbox{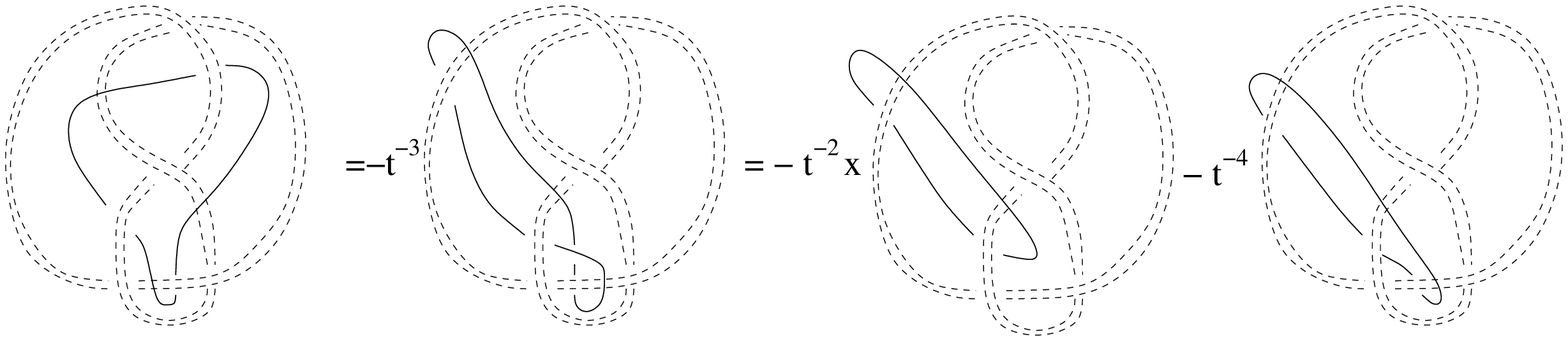}

Figure 4.5.    
\end{figure}

It follows that $B=x^2y+t^{-4}x^2z+(1-t^4)y+4t^{-2}x^2-t^{-2}-t^2$
After replacing $y^2$ and $z^2$ with their
 values given in Lemma \ref{yz} in A and
adding  A and B we obtain the desired formula.
\end{proof}

To simplify the computations we denote $Y=t^2y+1$, $Z=t^{-2}z+1$.
Also recall the polynomials $S_n(x)$, with $S_0(x)=1$, $S_1(x)=x$,
$S_{n+1}(x)=xS_n(x)-S_{n-1}(x)$, $n\in {\mathbb Z}$. 

\begin{prop}\label{(1,q)}
For any integer $q$, the following formula holds
\begin{eqnarray*}
t^{-q}(1,q)_T=(t^2S_{q+2}(x)-t^{-2}S_{q-2}(x))Y+(t^2S_q(x)-t^{-2}S_{q-4}(x))Z
\\+(t^2S_{q+2}(x)-t^{-2}S_{q-4}(x)).
\end{eqnarray*}
\end{prop}

\begin{proof}
For $q=0$ and $q=1$ this follows from Lemma \ref{(1,0)}. The product-to-sum
formula yields the recursive relation
\begin{eqnarray*}
(1,q+1)_T=tx(1,q)_T-t^2(1,q-1)_T
\end{eqnarray*}
and an inductive argument proves the formula in general.
\end{proof}

\begin{lemma}\label{YZaction}
The action of the skein algebra of the boundary torus on the
skein module of the figure-eight knot complement is determined  by
\begin{alignat*}{1}
(1,q)_T Y&=-t^2(1,q+2)_T+t^{q+2}T_{q+2}(x)Y+t^{q+2}T_q(x)\\
(1,q)_T Z&=-t^{-2}(1,q-2)_T+t^{q-2}T_{q-2}(x)Z+t^{q-2}T_{q}(x),
\end{alignat*}
where $q$ is any integer.
\end{lemma}

\begin{proof}
In the proof we prefer to work with $y$ instead of $Y=t^2y+1$.
It suffices to prove the formula for the action on $y$,
the action on $z$ is then determined by taking the mirror
image. 
The idea is to compute ``by hand'' $(1,0)_T y$ and $(1,1)_T y$,
and then use the product-to-sum formula to write the
recurrence:
\begin{eqnarray*}
(1,q+1)_T y=tx(1,q)_T y-t^2(1,q-1)_T y.
\end{eqnarray*}
We explain how $(1,0)_T y $ is determined and leave the case
of $(1,1)_T y$ to the reader. 
We begin by resolving the crossing in $(1,0)_Ty$ indicated
by an arrow  in Figure 4.6.
\begin{figure}[htbp]
\centering
\leavevmode
\epsfxsize=5in
\epsfysize=2.5in
\epsfbox{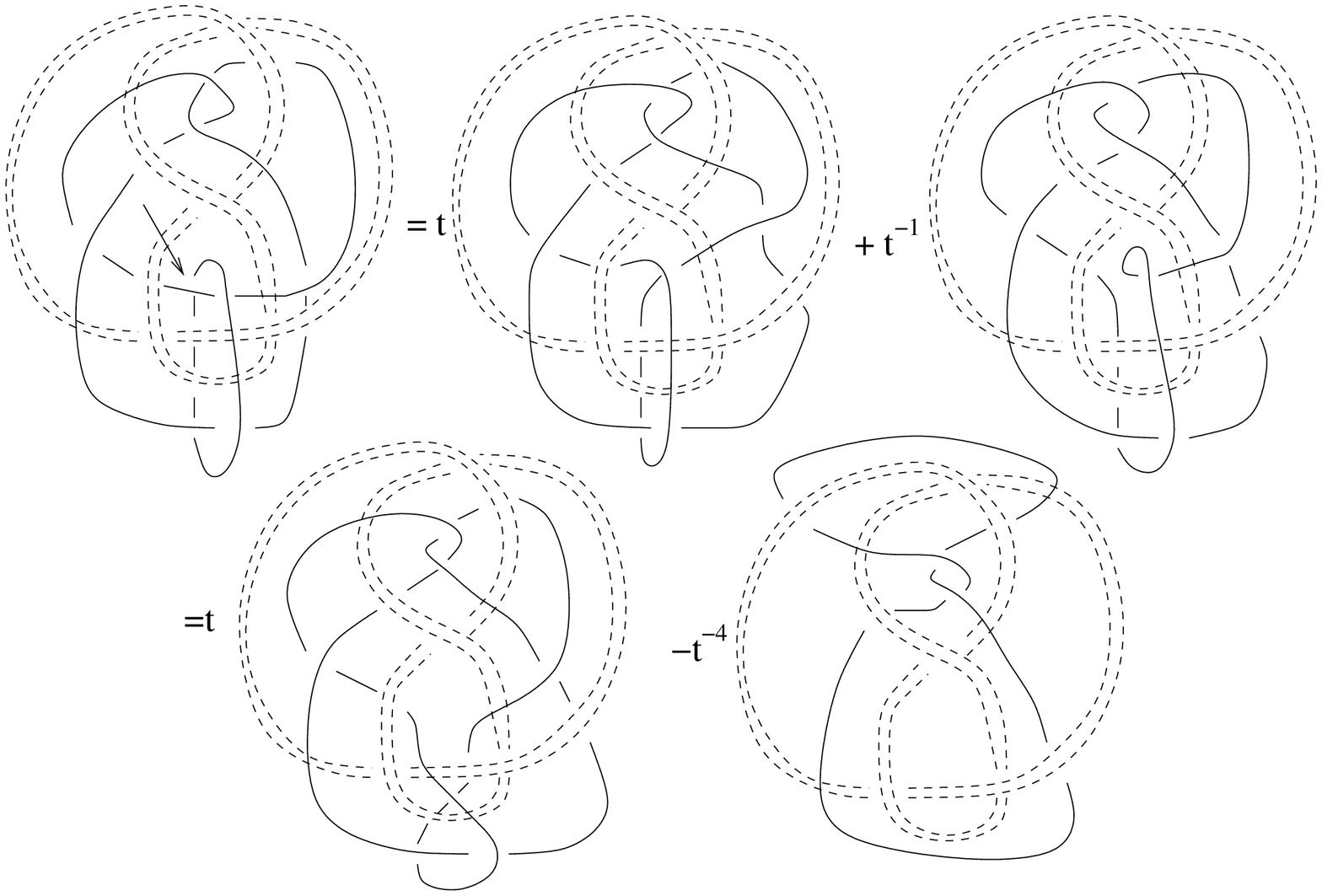}

Figure 4.6.    
\end{figure}
After an isotopy, the second term becomes $-t^{-4}t^6y=-t^2y$. The 
first term is nothing but $t^2(-t^{-3}(1,1)_T*(0,1)_T)+D$ where
$D$ is the first diagram from Figure 4.7. The remaining part
of  Figure 4.7 computes the value of $D$. To conclude, the last
two diagrams are computed by a method already presented in the
proof of Lemma \ref{(1,0)} (at the end of Figure 4.4). They are
found to be $-t^4y(-t^{-2}x^2-t^{-4}z)$ and respectively $-t^2(-t^{-2}x^2
-t^{-4}z)$. By adding all terms and using Lemma \ref{yz} and the product
to-sum formula we obtain
\begin{eqnarray*}
(1,0)_T y=
-(1,2)_T-t^{-2}(1,0)_T+t^2T_2(x)y+xT_1(x)
\end{eqnarray*}
Finally, replacing $t^2y+1$ by $Y$ we obtain the desired result. 

\begin{figure}[htbp]
\centering
\leavevmode
\epsfxsize=5in
\epsfysize=1.5in
\epsfbox{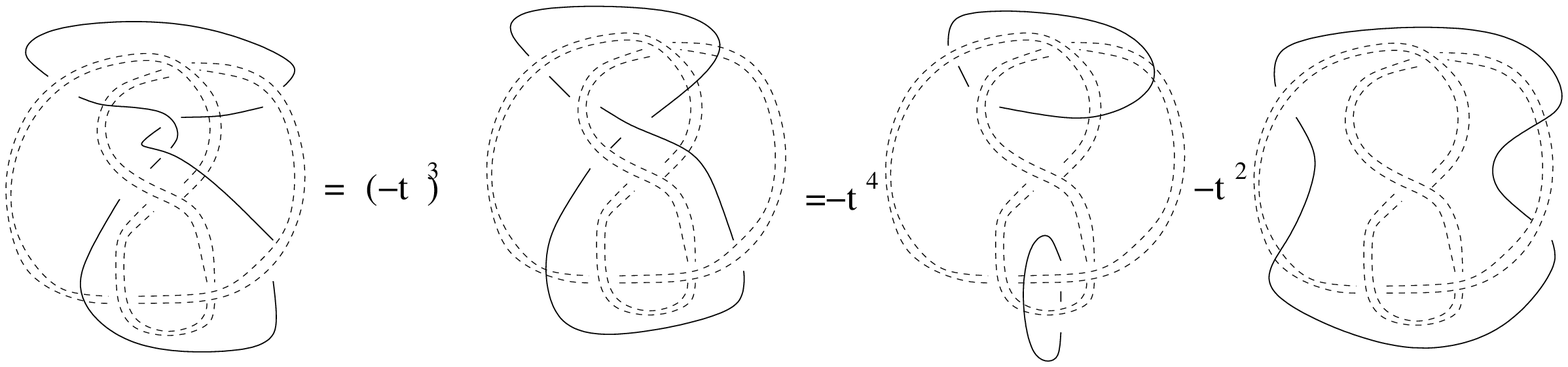}

Figure 4.7.    
\end{figure}
\end{proof}

\begin{remark}
 All  the information
about  the action of $K_t({\mathbb T}^2\times I)$ on
the skein module of the figure-eight knot complement is contained 
in Proposition \ref{(1,q)} and Lemma \ref{YZaction}.
The action of some 
$(p,q)_T\in K_t({\mathbb T}^2\times I)$
 on a basis element of the form $x^nY$ or $x^nZ$ can be computed
using the product-to-sum formula. 
\end{remark}

\section{The generators of the peripheral ideal}

We are now  done  with all necessary skein computations.
From this moment  on our work will be  based on 
  noncommutative trigonometry
 at the level of the Kauffman bracket skein
algebra of the torus (hence in the noncommutative torus).    
Here is an immediate consequence of Lemma \ref{(1,0)} that
will be important in our computations.
\begin{lemma}\label{t5t3t1}
The following formulas are true
\begin{alignat*}{1}
[S_5(x)-(1+t^4+t^{-4})S_1(x)]Y&=t^{-5}(1,3)_T-t^3(1,-1)_T+t^4(0,1)_T-
(0,5)_T\\
[S_5(x)-(1+t^4+t^{-4})S_1(x)]Z&=t^{5}(1,-3)_T-t^{-3}(1,1)_T+t^{-4}(0,1)_T-
(0,5)_T.
\end{alignat*}
\end{lemma}

We stress out that the   lemma states  that the left-hand sides are equal to
the images in the skein module of the knot complement of
the elements of the skein algebra of the torus  written here on
the right-hand side.

The next result, combined with Proposition \ref{(1,q)}, tells us
what the image of $(2,q)_T$ is in the skein module of the 
figure-eight knot complement. 

\begin{lemma}\label{2q}
For any integer $q$, the following identity holds
\begin{eqnarray*}
& & t^{-2q}(2,q)_T=[t^8(0,q+4)_T+(t^4+1)(0,q+2)_T+(0,q)_T]Y\\
& & \quad +[t^{-8}(0,q-4)_T+(t^{-4}+1)(0,q-2)_T+(0,q)_T]Z\\
& &  \quad +t^{-q}[-t^6(1,q+4)_T-
t^{-6}(1,q-4)_T]\\
& & \quad +(t^8+1)(0,q+2)_T+
(t^4+t^{-4}+1)(0,q)_T+(t^{-8}+1)(0,q-2)_T.
\end{eqnarray*}
\end{lemma}

\begin{proof}
By the product-to-sum formula 
\begin{eqnarray*}
t^{-2q}(2,q)_T+(0,q)_T=t^{-q}(1,q)_T*(1,0)_T.
\end{eqnarray*}
Using Lemma \ref{(1,0)} we see that this is further
equal to
\begin{eqnarray*}
& & t^{-q}(1,q)[(t^2(0,2)_T+(t^2+t^{-2}))Y\\
& & +(t^{-2}(0,2)_T+(t^2+t^{-2}))Z+
(t^2+t^{-2})(0,2)_T+(t^2+t^{-2})]
\end{eqnarray*}
or
\begin{eqnarray*}
& & =t^{-q}[t^4(1,q+2)_TY+(1,q-2)_TY+(t^2+t^{-2})(1,q)_TY\\
& & +(1,q+2)_TZ+t^{-4}(1,q-2)_TZ+(t^2+t^{-2})(1,q)_TZ\\
& & +(t^4+1)(1,q+2)_T+(1+t^{-4})(1,q-2)_T+(t^2+t^{-2})(1,q)_T].
\end{eqnarray*}
Applying Lemma \ref{YZaction} 
we obtain the desired formula.
\end{proof}

\begin{prop}\label{(2,q)}
For any integer $q$, one has
\begin{eqnarray*}
& & t^{-2q}(2,q)_T=[t^{-8}S_{q-4}(x)-t^{-4}S_q(x)+t^{-4}S_{q-2}(x)-
t^8S_{q+2}(x)]Y\\
& &\quad +[t^{-8}S_{q-4}(x)-t^{4}S_q(x)+t^{4}S_{q-2}(x)-t^8S_{q+2}(x)]Z\\
&& \quad
+t^{-q}[-t^6(1,q+4)_T-t^{-6}(1,q-4)_T+t^4(1,q+2)_T+t^{-4}(1,q-2)_T\\
& & \quad +(t^2+
t^{-2})(1,q)_T]
 +[-t^8T_{q+4}(x)-t^4T_{q+2}(x)-t^8S_q(x)\\
& & \quad +t^{-8}S_{q-2}(x)-t^{-4}T_{q-2}(x)
-t^{-8}T_{q-4}(x)].
\end{eqnarray*}
\end{prop}

\begin{proof}
Recall that $(0,q)_T=T_q(x)$ and $T_q(x)=S_q(x)-S_{q-2}(x)$ for any
integer $q$. 
The proposition then  follows from  Lemma \ref{2q} by  applying the
formulas from   Proposition \ref{(1,q)} to the terms
that contain $(1,q\pm 2)_T$ and $(1,q)_T$. 
\end{proof}

\begin{prop}\label{peri2}
The elements
\begin{eqnarray*}
& & t^{-6}(2,3)_T-t^6(2,-1)_T+t^3(1,7)_T-t(1,5)_T\\
& &\quad  +(-t^{11}+t^3-t^{-1}-t^{-5})
(1,3)_T
 +(t^{9}-t^5-t^{-7})(1,1)_T\\
& & \quad +(-t^{11}+2t^7+t^3-t^{-1}+t^{-9})(1,-1)_T+
(t^{13}+t)(1,-3)_T\\
& & \quad -t^{-1}(1,-5)_T
 +t^8(0,7)_T+(-2t^8+t^4-t^{-4})(0,5)_T\\
& & \quad +(-t^{12}+t^8-t^4-1+t^{-4})(0,3)_T
+(t^{12}-t^8+1+t^{-4})(0,1)_T
\end{eqnarray*}
and
\begin{eqnarray*}
& & t^{6}(2,-3)_T-t^{-6}(2,1)_T+t^{-3}(1,-7)_T-t^{-1}(1,-5)_T\\
& &\quad  +(-t^{-11}+t^{-3}-t^{}-t^{5})
(1,-3)_T
 +(t^{-9}-t^{-5}-t^{7})(1,-1)_T\\
& & \quad +(-t^{-11}+2t^{-7}+t^{-3}-t^{}+t^{9})(1,1)_T+
(t^{-13}+t^{-1})(1,3)_T\\
& & \quad -t(1,5)_T
 +t^{-8}(0,7)_T+(-2t^{-8}+t^{-4}-t^{4})(0,5)_T\\
& & \quad +(-t^{-12}+t^{-8}-t^{-4}-1+t^{4})(0,3)_T
+(t^{-12}-t^{-8}+1+t^{4})(0,1)_T
\end{eqnarray*}
are in the peripheral ideal of the figure-eight knot.
\end{prop}

\begin{proof}
Using Proposition \ref{(2,q)} we obtain
\begin{eqnarray*}
&& t^{-6}(2,3)_T=[t^{-8}S_{-1}(x)-t^{-4}S_3(x)+t^{-4}S_1(x)-t^{8}S_5(x)]Y\\
&& +[t^{-8}S_{-1}(x)-t^4S_3(x)+t^4S_1(x)-t^8S_5(x)]Z\\
&& +t^{-3}[-t^6(1,7)_T-t^{-6}(1,-1)_T+t^4(1,5)_T+t^{-4}(1,1)_T+
(t^2+t^{-2})(1,3)_T]\\
& & +[-t^8T_7(x)-t^4T_5(x)-t^8S_3(x)+t^{-8}S_1(x)-t^{-4}T_1(x)-
t^{-8}T_{-1}(x)].
\end{eqnarray*}
Using Lemma \ref{t5t3t1} and the fact that $S_{-n}=-S_{n-2}$,  $S_{-1}=0$, and
$T_{-n}=T_n$, we obtain
\begin{eqnarray*}
& & t^{-6}(2,3)_T=[-t^{-4}S_3(x)+(t^{-4}-t^{12}-t^4-t^8)S_1(x)]Y\\
& & \quad +[-t^{-4}S_3(x)+(-t^{12}-t^8)S_1(x)]Z\\
& & \quad +[-t^{3}(1,7)_T+t(1,5)_T+
(t^{-1}+t^{-5}-t^3)(1,3)_T+(t^{-7}+t^5)(1,1)_T\\
& & \quad +(t^{11}-t^{-9})(1,-1)_T-
t^{13}(1,-3)_T]+[-t^8(0,7)_T\\
& & \quad +(2t^8-t^4)(0,5)_T-t^8(0,3)_T-(t^8+t^4+t^{-4}+t^{12})(0,1)_T].
\end{eqnarray*}
A similar computation yields
\begin{eqnarray*}
& & t^6(2,-1)_T=[-t^{-4}S_3(x)-(t^{12}+1)S_1(x)]Y\\
& & \quad +[-t^{-4}S_3(x)-(t^{12}+t^8)S_1(x)]Z\\
& & \quad 
+[-t^{11}(1,3)_T-t^{-1}(1,-5)_T+t^{11}(1,1)_T+t(1,-3)_T+(t^7+t^3)(1,-1)_T]\\
& & \quad +[(-t^{12}-1)(0,3)_T-t^{-4}(0,5)_T+(t^{-4}-t^8)(0,1)_T].
\end{eqnarray*}
Now when we subtract the two, everything is in the image of the 
$K_t({\mathbb T}^2\times I)$ except for 
$(t^{-4}-t^4-t^8+1)S_1(x)Y+(t^{-4}-t^4)S_3(x)Z$. We can rewrite this
as
\begin{eqnarray*} 
(-t^6-t^{-2})[(t^2S_1(x)-t^{-2}S_{-3}(x))Y-(t^2S_{-1}(x)-t^{-2}S_{-5}(x))Z]
\end{eqnarray*}
and then we recognize the $Y,Z$-parts in the 
formula for $(1,-1)_T$ as given by Proposition 
\ref{(1,0)}. Using that formula and moving everything to the right
we obtain an element in $K_t({\mathbb T}^2\times I)$ that is mapped
to zero in the skein module of the knot complement. This
is exactly the first element from the statement. 

The second
element  is the  mirror image of the first, and since the figure eight knot is
amphichiral, it also belongs to the peripheral ideal.
\end{proof}

For a nonnegative integer $n$ let us denote by
${V}(n)$ the intersection of the peripheral
ideal $I_t(K)$ with the ${\mathbb C}(t)$-submodule
of $K_t({\mathbb T}^2\times I)$ that is generated by elements
of  the form $(p,q)_T$, $0\leq p\leq n$. We point out that
${V}(n)$  is not an ideal. 

\begin{lemma}\label{xYxZ}
The elements $xY$ and $xZ$ are in the image of  ${V}(2)$ in the 
the skein module of the figure-eight knot complement. 
\end{lemma}

\begin{proof}
If in the formula for $(2,-1)_T$ obtained in the proof of 
Proposition \ref{peri2} we add on both
sides $t^{-3}(1,-1)_T+t^{-7}(1,1)_T$ and  expand these two
terms on the right-hand side using Lemma \ref{(1,q)} we obtain
\begin{eqnarray*}
&&t^{6}(2,-1)_T+t^{-3}(1,-1)_T+t^{-7}(1,1)_T=
-t^{-4}(t^{16}-1)S_1(x)Y\\
& & -t^{-4}(t^{16}-1)(t^{-4}-1)S_1(x)Z-
t^{-1}(1,-5)_T+t^{-1}(1,-3)_T
+(t^7+t^3)(1,-1)_T\\
&  & +t^9(1,1)_T-t^{11}(1,3)_T -
t^{-4}(0,5)_T+(t^{-4}-t^{12}-1)(0,3)_T\\
& & +(-t^{8}-1-t^{-4}-t^{-8})(0,1)_T.
\end{eqnarray*}
Similarly
\begin{eqnarray*}
&&t^{-6}(2,1)_T+t^{3}(1,1)_T+t^{7}(1,-1)_T=
t^{-12}(t^{16}-1)S_1(x)Z\\
& & t^{-12}(t^{16}-1)(t^{4}-1)S_1(x)Y-
t^{}(1,5)_T+t(1,3)_T\\
& & +(t^{-7}+t^{-3})(1,1)_T+t^{-9}(1,-1)_T-t^{-11}(1,-3)_T\\
& & -
t^{4}(0,5)_T+(t^{4}-t^{-12}-1)(0,3)_T+(-t^{-8}-1-t^{4}-t^{8})(0,1)_T.
\end{eqnarray*}
Recall that $S_1(x)=x$. Multiplying the second equation by $t^4-t^8$,
adding the two,
 and dividing by $(t^{16}-1)(t^{-8}-t^{-4}+1)$  we obtain
\begin{eqnarray*}
& & xY=\frac{1}{(t^{16}-1)(1-
t^{-4}+t^{-8})}[(t^{-2}-t^{2})(2,1)_T+t^{6}(2,-1)_T+
(t^5-t^9)(1,5)_T\\
& & \quad  +(t^{11}+t^9-t^5)(1,3)_T
+(-t^{11}-t^9+t^7-t^5-t^{-3}+t^{-7})(1,1)_T\\
& & \quad +(t^{-3}-t^{15}+t^{11}-t^{7}-t^3+
t^{-1}-t^{-5})(1,-1)_T\\
& & \quad +(-t^{-1}-t^{-3}+t^{-7})(1,-3)_T
+t^{-1}(1,-5)_T+(-t^{12}+t^8+t^{-4})(0,5)_T\\
& & \quad 
+(2t^{12}-2t^8+t^4-2t^{-4}+t^{-8}+1)(0,3)_T+(-t^{16}+t^8+t^4+2t^{-4}+t^{-8})(0,1)_T].
\end{eqnarray*}
Taking the mirror image we obtain
\begin{eqnarray*}
& & xZ=\frac{1}{(t^{-16}-1)(1-
t^{4}+t^{8})}[(t^2-t^{-2})(2,-1)_T+t^{-6}(2,1)_T
+(t^{-5}-t^{-9})(1,-5)_T\\
& & \quad +(t^{-11}+t^{-9}-t^{-5})(1,-3)_T
+(-t^{-11}-t^{-9}+t^{-7}+t^{-5}-t^{3}+t^{7})(1,-1)_T
\\ & & \quad +
(t^3-t^{-15}+t^{-11}-t^{-7}-t^{-3}+
t-t^{5})(1,1)_T\\
& & \quad 
+(-t-t^{3}+t^{7})(1,3)_T+t(1,5)_T+(t^{-12}+t^{-8}+t^{4})(0,5)_T\\
& & \quad +(2t^{-12}-2t^{-8}+t^{-4}-2t^4+t^8+1)(0,3)_T
+(-t^{-16}+t^{-8}+t^{-4}+
2t^{4}+t^{8})(0,1)_T].
\end{eqnarray*}
\end{proof}

From now on $xY$ and $xZ$ will be identified with the above preimages
in the skein algebra of the boundary torus. 

\begin{prop}\label{peri31}
The element 
\begin{eqnarray*}
& & -t(3,1)_T+(t^{16}-1)(t^4+1)t^{-2}(0,2)_T(xY)
+(t^{16}-1)(t^4+1)(t^{-6}-t^{-2})(xY)\\
& & +(t^{16}-1)t^{-2}(0,2)_T(xZ)
+(t^{16}-1)(t^{-6}-t^{-10})(xZ)\\
& & -t^{12}(2,5)_T
+t^8(2,3)_T+(1+t^4-t^{10})(2,1)_T+(t^{-4}-t^{12})(2,-1)_T-t^{-8}(2,-3)_T\\
& & -t^{15}(1,5)_T+(t^{-3}-t^9-t^{17})(1,3)_T+(t^{-5}-t^7-2t^{11}-2t^{15}+
t^{-1})
(1,1)_T\\
& & (t^{-7}-t^{-3}-t^{-1}-t^5-t^{13})(1,-1)_T+
(t^{-9}-t^{-5}-t^3)(1,-3)_T-t^5(1,-5)_T\\
& & -(t^2-t^{-2})(0,3)_T
+(t^{18}+t^{14}+2t^{10}-t^8+t^6-2t^{-6}-t^{-10})(0,1)_T
\end{eqnarray*}
is in the peripheral ideal of the figure-eight knot.
\end{prop}

Note that $xY$ and $xZ$ were shown to be in the image
of the skein algebra of the boundary, and it follows that
$(0,2)_T(xY)$ and $(0,2)_T(xZ)$ are in the image of this algebra, too.

\begin{proof}
By the product-to-sum formula we have
\begin{eqnarray*}
t(3,1)_T-(1,0)_T*(2,1)_T+t^{-1}(1,1)_T=0.
\end{eqnarray*}
Substitute in this $(1,0)_T$ using Lemma \ref{(1,0)}, and
$(2,1)_T$ using the formula deduced in  the proof of
Proposition \ref{xYxZ}.   
\end{proof}

\begin{lemma}\label{02Y}
The following equalities hold
\begin{alignat*}{1}
(0,2)_TY&=(0,1)_T(xY)-2Y\\
(0,2)_TZ&=(0,1)_T(xZ)-2Z.
\end{alignat*}
\end{lemma}

\begin{proof}
Since $xY=(0,1)_TY$ and $(0,0)_T=2$ the conclusion follows
from 
\begin{eqnarray*}
(0,2)_TY=(0,1)_T*(0,1)_TY-(0,0)_T
\end{eqnarray*}
and its mirror image.
\end{proof}

\begin{lemma}\label{YZ}
One has the following formulas
\begin{eqnarray*}
& & Y=\frac{1}{2(t^{-8}-t^4)(t^4-1)}[(t^4+t^{-4}+1-2t^8-t^{-8})(0,1)_T(xY)\\
& & \quad +(t^{-4}+1-t^4-t^8)(0,1)_T(xZ)-(2,0)_T-t^6(1,4)_T-t^{-6}(1,-4)_T\\
& & \quad +t^4(1,2)_T+t^{-4}(1,-2)_T+(t^6-t^{-6})(1,0)_T-(t^8+t^{-8})(0,4)_T\\
& & \quad +(t^{-4}-t^{4}+2-t^{8}+t^{-8})(0,2)_T+(2t^{-4}+2-2t^{8})]
\end{eqnarray*}
and
\begin{eqnarray*}
& & Z=\frac{1}{2(t^{8}-t^{-4})(t^{-4}-1)}[(t^{-4}+t^4+1-2t^{-8}-t^8)
(0,1)_T(xZ)\\
& & \quad 
+(t^{4}+1-t^{-4}-t^{-8})(0,1)_T(xY)-(2,0)_T-t^6(1,4)_T-t^{-6}(1,-4)_T\\
& & \quad 
+t^4(1,2)_T+t^{-4}(1,-2)_T+(t^{-6}-t^{6})(1,0)_T-(t^8+t^{-8})(0,4)_T\\
& & \quad +(t^4-t^{-4}+2-t^{-8}+t^8)(0,2)_T+(2t^4+2-2t^{-8})]
\end{eqnarray*}
\end{lemma}

\begin{proof}
Using the product-to-sum formula we obtain
\begin{eqnarray*}
& & (2,0)_T=-(t^8+t^{-8})(0,2)_TY-(2t^{-4}+t^8+t^{-8})Y-(t^8+t^{_8})(0,2)_TZ\\
& & -(2t^4+t^8+t^{-8})Z-t^6(1,4)_T-t^{-6}(1,-4)_T+t^4(1,2)_T+t^{-4}(1,-2)_T
\\
& & +(t^2+t^{-2})(1,0)_T
 -t^8(0,4)_T-t^4(0,2)_T-t^8-t^{-8}-t^{-4}(0,2)_T-t^{-8}(0,4)_T\\
& & \quad \quad =-(t^6+t^2+t^{-2}-t^{-6})(t^{-2}-t^2)Y-(t^{-6}+t^{-2}+t^2-t^6)(t^2+t^{-2})Z
\\
& & -(t^8+t^{-8})(0,1)_T(xY)-(t^8+t^{-8})(0,1)_T(xZ)\\
& & -t^6(1,4)_T-t^{-6}(1,-4)_T+t^4(1,2)_T+t^{-4}(1,-2)_T
+(t^2+t^{-2})(1,0)_T\\
& & -(t^8+t^{-8})(0,4)_T-(t^4+t^{-4})(0,2)_T-(t^8+t^{-8}).
\end{eqnarray*}
If we combine this with
\begin{alignat*}{1}
(1,0)_T=&t^2(0,1)_T(xY)+t^{-2}(0,1)_T(xZ)+(t^{-2}-t^2)Y+(t^2-t^{-2})Z\\
 &+(t^2+t^{-2})(0,2)_T+(t^2+t^{-2}),
\end{alignat*}
we obtain a system of two equations in the unknowns $Y$ and $Z$. 
Solving it we obtain the formulas from the statement.
\end{proof}

\begin{cor}
The Kauffman bracket skein module of the complement of the figure-eight
 knot is peripheral. Moreover, it is the image of ${V}(2)$.
\end{cor}
\begin{proof}
$Y$ and $Z$ are in ${V}(2)$ by Lemma \ref{YZ}. Since 
${V}(2)$ is closed under multiplication on the 
left by $(0,1)_T$ it follows that $x^n$, $x^nY$ and $x^nZ$ are in the
image of ${V}(2)$ for all nonnegative integers $n$. But
by \cite{bullof} these elements form a basis of the skein module 
of the knot complement.
\end{proof}

\begin{prop}\label{peri32}
The element
\begin{eqnarray*}
& & -t^3(3,0)_T-t^{-7}(1,4)_TY-t^{-13}(1,-2)_TY-(t^7+t^{-9}+t^{-5})(1,2)_TY\\
& & -
(t^5+t^{-11}+t^{-7})(1,0)_TY
-t^{-7}(1,4)_TZ-t^{-13}(1,-2)_TZ\\
& & -(t^7+t^{-9}+t^3)(1,2)_TZ
-
(t^5+t^{-11}+t)(1,0)_TZ\\
& & 
 -t^7(2,4)_T-t^{-13}(2,-4)_T+t^3(2,2)_T+t^{9}(2,-2)_T+(t^{-1}+t^{-5})(2,0)_T
\\
& & -(t^9-t^{-3})(1,4)_T-(t^3+t^{-3}+t^{-9})(1,-2)_T-t^{-5}(1,6)_T-
t^{-15}(1,-4)_T\\
& & +(t^{-9}-t^3)(1,2)_T
+(t^{-11}-t)(1,0)_T-t^{-1}(0,6)_T+t^{-1}(0,4)_T+t^{-1}(0,2)_T+2t^3.
\end{eqnarray*}
is in the peripheral ideal of the figure-eight knot.
\end{prop}


\begin{proof}
The proof is analogous to the one given to Proposition \ref{peri31}; 
in this case one  starts with
\begin{eqnarray*}
(3,0)_T-(1,0)_T*(2,0)_T+(1,0)_T=0.
\end{eqnarray*}
\end{proof}



\begin{theorem}
The four elements described in Propositions \ref{peri2}, \ref{peri31},
and \ref{peri32} generate the peripheral ideal of the figure-eight knot.
\end{theorem}

\begin{proof}
(1) First let us show that there are no nontrivial elements in
${V}(1)\cap I_t(K)$. Assume the contrary, namely that
there exist finitely many $c_q$, not all equal to zero,
such that $\sum_qc_qt^{-q}(1,q)_T$ is in the peripheral ideal.
The fact that $S_n(x)$, $S_n(x)Y$, $S_n(x)Z$, $n\geq 0$ form a basis
of the skein module of the knot complement, combined with Lemma \ref{(1,q)} 
implies that 
\begin{alignat*}{1}
\sum_qc_q(t^2S_{q+2}(x)-t^{-2}S_{q-2}(x))&=0\\
\sum_qc_q(t^2S_{q}(x)-t^{-2}S_{q-4}(x))&=0\\
\sum_qc_q(t^2S_{q+2}(x)-t^{-2}S_{q-4}(x))&=0.
\end{alignat*}
as the first is the coefficient of $Y$, the second is the coefficient
of $Z$ and  the third is the free term in the formula of $\sum_qc_q
t^{-q}(1,q)_T$.
Subtracting the first from the last and the second from the last we
obtain
\begin{eqnarray*}
t^{-2}\sum_qc_qT_{q-2}(x)=0 \quad \mbox{and} \quad 
t^2\sum_qc_qT_{q+2}(x)=0.
\end{eqnarray*}
Since the Chebyshev polynomials $T_n, n\geq 0$ form a basis of the 
vector   space of polynomials, it follows from the 
first relation that $c_{q+2}=c_{-q+2}$, and from the second that  
$c_{q-2}=c_{-q-2}$, for all $q$. This is impossible unless all
$c_q$ are equal to zero, a contradiction. Hence the claim is true.

(2) Now let us show that 
\begin{eqnarray*}
{V}(2)\cap I_t(K)=\{\alpha g_1+\beta g_2, \: \alpha,\beta \in
{\mathbb C}(t)[x]\}
\end{eqnarray*}
where $g_1$ and $g_2$ are the two elements exhibited in the
statement of Proposition \ref{peri2}.

Reduce the  ${\mathbb C}(t)[x]$-modules modulo ${V}(1)$. Then
\begin{eqnarray*}
t^{-6}(2,3)_T-t^6(2,-1)_T\equiv
0 \quad \mbox{and} \quad t^6(2,-3)_T-t^{-6}(2,1)_T\equiv 0.
\end{eqnarray*}
Multiplying on the left by $(0,1)_T$  and using the product-to-sum
formula we can prove inductively that for all integers $p$,
\begin{eqnarray*}
(2,2p+1)_T-t^{12}(2,2p-3)_T\equiv 0 \quad (\mbox{mod }{V}(1))
\end{eqnarray*}
and
\begin{eqnarray*}
(2,2p+2)_T+t^4(2,2p)_T-t^{12}(2,2p-2)_T-t^{6}(2,2p-4)_T\equiv 0\quad
 (\mbox{mod } {V}(1)).
\end{eqnarray*}
Assume that there exists an element in ${V}(2)\cap I_t(K)$ not
generated by $g_1$ and $g_2$. Subtracting terms generated by
$g_1$ and $g_2$ of the form
$(2,2p+1)_T-t^{12}(2,2p-3)_T+v_{2p+1}$ or $
(2,2p+2)_T+t^4(2,2p)_T-t^{12}(2,2p-2)_T-t^{6}(2,2p-4)_T+v_{2p+2}$
with $v_{2p+1},v_{2p+2}\in {V}(1)$, we can reduce
such  an element to one that is a linear combination
of $(2,2)_T,(2,1)_T,(2,0)_T, (2,-1)_T, (2,-2)_T$.
Examining the formulas in Propositions \ref{(1,q)} and \ref{(2,q)}
we see that parity implies that the ``even'' part and
the ``odd'' part can be separated. This means  that there must be an element in
${V}(2)\cap I_t(K)(\mbox{mod }{V}(1))$ of the form $a(2,1)_T+b(2,-1)_T$, and
one of the form $c(2,2)_T+d(2,0)_T+e(2,-2)_T$. We will show that $a=b=c=d=e=0$.

By Proposition \ref{(2,q)}, the fact that $a(2,1)_T+b(2,-1)_T$
is in the peripheral ideal modulo $V(1)$  means that
\begin{eqnarray*}
& & [-(t^{-4}a+t^4b)S_3(x)-(t^{12}a+a+t^{-12}b+t^{-8}b)S_1(x)]Y\\
& & \quad +[-(t^{-4}a+t^4b)S_3(x)-(t^{12}a+t^8a-t^{-12}b+b)S_1(x)]Z
\end{eqnarray*}
is in ${V}(1)$. Let it be equal to $\sum_q c_qt^{-q}(1,q)_T$. 
By parity we see that we can assume $c_q=0$ for $q$ even.
Since 
\begin{eqnarray*}
& & \sum_q c_qt^{-q}(1,q)_T=(t^2c_qS_{q+2}(x)-t^{-2}c_qS_{q-2}(x))Y\\
& & +(t^2c_qS_q(x)-t^{-2}c_qS_{q-4}(x))Z
+c_q(t^2S_{q+2}(x)-t^{-2}S_{q-4}(x)).
\end{eqnarray*}
$q$ can only be equal to $3,1,-1,-3$. Indeed, if we denote by $q_0$
and $q_1$ the largest and the smallest of the $q$'s, then
for $S_{q_0+2}(x)Y$ and $S_{q_1-2}(x)Y$ to cancel out one should have
$q_0+2=-(q_1-2)-2$, and for $S_{q_0}(x)Z$ and $S_{q_1-4}(x)Z$
to cancel out one should have $q_0=-(q_1-4)-2$. So on the one
hand $q_0=-q_1-2$ and on the other $q_0=-q_1+2$, and these two
cannot hold simultaneously.

Let us now look at 
\begin{eqnarray*}
c_3t^{-3}(1,3)_T+c_1t^{-1}(1,1)_T+c_{-1}t(1,-1)_T+
c_{-3}t^{3}(1,-3)_T.
\end{eqnarray*}
For the coefficients of $Y$ and $Z$ to be polynomials of third degree
in $x$ one should have $c_3=c_{-3}=0$.  So we are left with
$c_1t^{-1}(1,1)_T+c_{-1}t(1,-1)_T$. Equating the coefficients of
the basis elements $S_3(x)Y$, $S_1(x)Y$, $S_3(x)Z$ and $S_1(x)Z$ 
we obtain
\begin{eqnarray*}
& & c_1t^2=-(t^{-4}a+t^4b)\\
& & (-c_1t^{-2}+c_{-1}t^2+c_{-1}t^{-2})=-(t^{12}a+a+t^{-12}b+t^{-8}b)\\
& & c_{-1}t^{-2}=-(t^{-4}a+t^4b)\\
& & (t^2c_{-1}+t^{-2}c_{-1})=-(t^{12}a+t^8a-t^{-12}b+b).
\end{eqnarray*}
This is a linear homogeneous system in the unknowns $a,b,c_1,c_{-1}$.
An easy computer check shows that the determinant is 
 different from  zero, so the system has the unique solution
$a=b=c_1=c_{-1}=0$.

The argument for
$c(2,2)_T+d(2,0)_T+e(2,-2)_T$ is similar, only more laborious, and
we prefer to skip it.
This proves that  ${V}(2)\cap I_t(K)$ contains no other elements than
those generated by $g_1$ and $g_2$. 

(3) In Propositions \ref{peri31} and \ref{peri32} we have
found elements in the peripheral ideal that are of the form
$(3,0)_T+v_0$ and $(3,1)_T+v_1$, with $v_0,v_1\in {V}(2)$.
The product-to-sum formula
\begin{eqnarray*}
(3,q+1)_T=t^3(0,1)_T*(3,q)_T-t^6(3,q-1)_T
\end{eqnarray*} 
allows us to construct, for any integer $q$, an element in the peripheral ideal
of the form $(3,q)_T+v_{3,q}$ with $v_{3,q}\in {V}(2)$. 
Also, the product-to-sum
formula 
\begin{eqnarray*}
(p+1,q)_T=t^{3q}(p-2,q)_T*(3,0)_T-t^{6q}(p-5,q)_T
\end{eqnarray*}
allows us to construct, for any integers $p\geq 4$ and  $q,$ an
element in the peripheral ideal of the form $(p,q)_T+v_{p,q}$, with
$v_{p,q}\in {V}(p-1)$. 
If there existed an element in the peripheral ideal that is
not generated  by the four exhibited in the statement, then
we could  reduce it using elements of the form $(p,q)_T+v_{p,q}$
to one in ${V}(2)$. But all elements in ${V}(2)\cap I_t(K)$ are
generated by $g_1$ and $g_2$. Hence the conclusion. 
\end{proof}

\section{The noncommutative  A-ideal of the figure-eight knot}

The noncommutative A-ideal is obtained from the peripheral
ideal by extending it to the ring of trigonometric polynomials
in the noncommutative torus, and then contracting it to
the quantum plane. Specifically, in the four generators exhibited
in Propositions \ref{peri2}, \ref{peri31}, and \ref{peri32} we
replace each $(p,q)_T$ by $t^{-pq}(l^pm^q+l^{-p}m^{-q})$
to obtain a Laurent polynomial in the noncommuting variables
$l$ and $m$ with coefficients in ${\mathbb C}(t)$. Then,
we multiply these Laurent polynomials on the left by appropriate
powers of $l$ and $m$ to obtain polynomials in ${\mathbb C}_t[l,m]$.
These four polynomials are the generators of the 
noncommutative A-ideal. 

In Lemma \ref{xYxZ} we saw that $xY$ and $xZ$ are the images
of two elements from the Kauffman bracket skein module of the torus.
We identified $xY$ and $xZ$ with their preimages. Now extend those
to the noncommutative torus, to obtain
\begin{eqnarray*}
& & xY=\frac{1}{(t^{16}-1)(1-t^{-4}+t^{-8})}[(-1+t^{-4})(l^2m+l^{-2}m^{-1})+
t^{8}(l^2m^{-1}+l^{-2}m)\\
& & \quad +(-t^4+1)(lm^5+l^{-1}m^{-5})
 +(t^8+t^6-t^2)(lm^3+l^{-1}m^{-3})\\
& & \quad +(-t^{10}-t^{8}+t^6-t^{-4}+t^{-8})(
lm+l^{-1}m^{-1})\\
& & \quad +(t^{-2}-t^{16}+t^{12}-t^{8}-t^4+1-t^{-4})(lm^{-1}-l^{-1}m)\\
& & \quad +(-t^2-1+t^{-4})(lm^{-3}+l^{-1}m^3) +t^4(lm^{-5}+l^{-1}m^5)\\
& & \quad + (-t^12+t^8+t^{-4})(m^5+m^{-5}) +(2t^{12}-2t^8+t^4-2t^{-4}+t^{-8}+1)
(m^3+m^{-3})\\
& & \quad +(-t^{16}+t^8+t^4+2+2t^{-4}+t^{-8})
(m+m^{-1})].
\end{eqnarray*}

and 

\begin{eqnarray*}
& & xZ=\frac{1}{(t^{-16}-1)(1-t^4+t^8)}[(t^4-1)(l^2m^{-1}+l^{-2}m)+ t^{-8}(l^2m+l^{-2}m^{-1})\\
& & \quad +(-t^{-4}+1)(lm^{-5}+l^{-1}m^5)+ (t^{-8}+t^{-6}-t^{-2})(lm^{-3}+l^{-1}m^3)\\
& & \quad +(-t^{-10}-t^{-8}+t^{-6}-t^{-4}-t^{4}+t^8)(lm^{-1}+l^{-1}m)\\
& & \quad +(t^2-t^{-16}+t^{-12}-t^{-8}-t^{-4}+1-t^4)(lm+l^{-1}m^{-1})\\
& & \quad +(-t^{-2}-1+t^4)(lm^3+l^{-1}m^{-3})+t^{-4}(lm^5+l^{-1}m^{-5})\\
& & +(t^{-12}+t^{-8}+t^4)(m^5+m^{-5}) +(2t^{-12}-2t^{-8}+t^{-4}-2t^4+t^8+1)(m^3+m^{-3})\\
& & \quad +(-t^{-16}+t^{-8}+t^{-4}+2+2t^{4}+t^{8})
(m+m^{-1})].
\end{eqnarray*}

Then in view of Lemma \ref{YZ} we can write
\begin{eqnarray*}
& & Y=\frac{1}{2(t^{-8}-t^4)(t^4-1)}[(t^4+t^{-4}+1-2t^8-t^{-8})(m+m^{-1})(xY)\\
& & \quad + (t^{-4}+1-t^4-t^8)(m+m^{-1})(xZ)-(l^2+l^{-2})\\
& & \quad -t^2(lm^4+l^{-1}m^{-4})-t^{-2}(lm^{-4}+l^{-1}m^4)\\
& & \quad +t^{2}(lm^2+l^{-1}m^{-2})+t^{-2}(lm^{-2}+l^{-1}m^2)\\
& & \quad +(t^6-t^{-6})(l+l^{-1})-(t^8+t^{-8})(m^4+m^{-4})\\
& & \quad +(t^{-4}-t^4+2-t^8+t^{-8})(m^2+m^{-2})+(2t^{-4}+2-2t^8)]
\end{eqnarray*}
and
\begin{eqnarray*}
& & Z=\frac{1}{2(t^{8}-t^{-4})(t^{-4}-1)}[(t^{-4}+t^{4}+1-2t^{-8}-t^{8})(m+m^{-1})(xZ)\\
& & \quad + (t^{4}+1-t^{-4}-t^{-8})(m+m^{-1})(xY)-(l^2+l^{-2})\\
& & \quad -t^{2}(lm^{4}+l^{-1}m^{-4})-t^{-2}(lm^{-4}+l^{-1}m^4)\\
& & \quad +t^{2}(lm^2+l^{-1}m^{-2})+t^{-2}(lm^{-2}+l^{-1}m^2)\\
& & \quad +(t^{-6}-t^{6})(l+l^{-1})-(t^8+t^{-8})(m^4+m^{-4})\\
& & \quad +(t^{4}-t^{-4}+2-t^{-8}+t^{8})(m^2+m^{-2})+(2t^{4}+2-2t^{-8})]
\end{eqnarray*}

From these formulas and  Theorem 5.11 we obtain the following 

\begin{theorem}
The noncommutative A-ideal of the figure-eight knot is generated
by the following four elements
\begin{eqnarray*}
& & t^{-18}l^3m^{14}+t^{8}l^2m^{14}-t^{-18}l^3m^{12}+(t^4+t^{-4}-2t^8)l^2m^{12}
 -t^{18}lm^{12}+t^{-40}l^4m^{10}\\
& & \quad +(-t^{-22}-t^{-18}+t^{-14}-t^{-6})
l^3m^{10}+ (t^8-t^4+t^{-4}-1-t^{12})l^2m^{10}\\
& & \quad +(t^{30}+t^{18})lm^{10}+ (-t^{-22}+t^{-6}-t^{-10})l^3m^8
+(-t^8+1+t^{-4}+t^{12})l^2m^8\\
& & \quad +(t^6+2t^{22}-t^{14}-t^{6}+t^{18})lm^8+t^{-20}l^4m^6\\
& & \quad +(-t^{-2}+t^{-22}-t^{-14}+2t^{-6}+t^{-10}+t^{-4})l^3m^6
+(1+t^{12}-t^{8}+t^{-4})l^2m^6\\
& & \quad +(t^{22}-t^{6}-t^{18}+t^{-4})lm^6
+(t^2+t^{-10})l^3m^4+(t^{-2}+t^8-t^4+t^{-4}-t^{12})l^2m^4\\
& & \quad +(-t^6-t^{10}+t^{14}-t^{22})lm^4+t^{16}m^4-t^{-10}l^3m^2
+(t^4-t^{-4}-2t^8)l^2m^2-t^{10}lm^2\\
& & \quad +t^8l^2+t^{10}l,\\
\end{eqnarray*}
\begin{eqnarray*}
& & t^{-8}l^2m^{14}+t^{18}lm^{14}-t^{-18}l^3m^{12}
 +(t^4+t^{-4}-2t^{-8})l^2m^{12}-t^{18}lm^{12}\\
& & \quad +(t^{-30}+t^{-18})l^3m^{10}+(t^4+t^{-8}-1-t^{-4}-t^{-12})l^2m^{10}\\
& & \quad +(t^{14}-t^{-22}-t^6-t^{18})lm^{10}+t^{40}m^{10} -t^{-36}l^4m^8\\
& & \quad +(-t^{-14}-t^{-26}+t^6+t^{-18}+2t^{-22})l^3m^8+(t^{-12}+t^4+1-t^{-8})l^2m^8\\
& & \quad +(-t^{22}-t^{10}+t^6)lm^8+(t^{-22}-t^{-18}-t^{-6})l^3m^6
+(1+t^{4}+t^{-12}-t^{-8})l^2m^6\\
& &  \quad +(-t^{24}+2t^6-t^2+t^{10})lm^6-t^{-20}m^6+t^{-16}l^4m^4\\
& & \quad+(-t^{-22}-t^{-6}-t^{-10}+t^{-14})l^3m^4+(t^4-t^{-12}+t^{-8}-t^{-4}-1)l^2m^4 \\
& & \quad +(t^{-2}+t^{10})lm^4-t^{-10}l^3m^2+(t^{-4}-2t^{-8}-t^{4})l^2m^2-t^{10}lm^2+t^{-8}l^2+t^{-10}l^3,
\end{eqnarray*}
\begin{eqnarray*}
& & l^3m^7[-t^{-2}(l^3m+l^{-3}m^{-1})+(t^{16}-1)(t^4+1)t^{-2}(m^2+m^{-2}
+t^{-4}-1)
(xY)\\
& & +(t^{16}-1)t^{-2}(m^2+m^{-2}+t^{-4}-t^{-8})(xZ)
 -t^{2}(l^2m^5+l^{2}m^{-5})
\\
& & +t^{2}(l^2m^3+l^{-2}m^{-3})
 +(t^{-2}+t^2-t^8)(l^2m+l^{-2}m^{-1})\\
& & +(t^{-2}-t^{14}(l^2m^{-1}+l^{-1}m^2)
-t^{-2}(l^2m^{-3}+l^{-2}m^3)\\
& &  -t^{10}(lm^5+l^{-1}m^{-5})+(t^{-6}-t^{6}-t^{14})(lm^3+l^{-1}m^{-3})\\
& & +(t^{-6}-t^6-2t^{11}-2t^{15}+t^{-1})(lm+l^{-1}m^{-1}) \\
& & + (t^{-6}-t^{-2}-1-t^{6}-t^{14})(lm^{-1}-l^{-1}m)
+(t^{-6}-t^{-2}-t^6)(lm^{-3}-l^{-1}m^3)\\
& & -t^{10}(lm^{-5}+l^{-1}m^5)
 -(t^2-t^{-2}(m^3-m^{-3})\\ & & +(t^{18}+
t^{14}+2t^{10}-t^8+t^6-2t^{-6}-t^{-10})(m+m^{-1})],
\end{eqnarray*}
and
\begin{eqnarray*}
& & 
l^3m^{10}[-t^3(l^3+l^{-3})+(-t^{-11}(lm^4+l^{-1}m^{-4})-t^{-11}(lm^{-2}+l^{-1}
m^2)\\
& & -(t^5+t^{-11}+t^{-7})(lm^2+l^{-1}m^{-2})
-(t^5+t^{-11}+t^{-7})(l+l^{-1}))Y\\
& & + (-t^{-11}(lm^4+l^{-1}m^{-4})-t^{-11}(lm^{-2}+l^{-1}m^2)\\
& & -(t^5+t^{-11}+t)
(lm^2+l^{-1}m^{-2})-(t^5+t^{-11}+t)(l+l^{-1}))Z\\
& & -t^{-1}(l^2m^4+l^{-2}m^{-4})-t^{-5}(lm^{-4}+l^{-1}m^4)+
t^{-1}(l^2m^2+l^{-2}m^{-2})\\
& & +(t^{-1}+t^{-5})(l^2+l^{-2})-(t^{5}-t^{-7})(lm^4+l^{-1}m^{-4})\\
& &- (t^5+t^{-1}+t^{-7}(lm^{-2}+l^{-1}m^2)-t^{-11}(lm^6+l^{-1}m^{-6})\\
& & -t^{-11}(lm^{-4}+l^{-1}m^4)+(t^{-11}-t^{-1})(lm^2+l^{-1}m^{-2})\\
& & +(t^{-11}-t)(l+l^{-1})-t^{-1}(m^6+m^{-6})+
t^{-1}(m^4+m^{-4})+t^{-1}(m^2+m^{-2})+2t^3].
\end{eqnarray*}
\end{theorem}

\section{The A-ideal determines the Jones polynomial}

For a knot $K$ we denote by $\kappa _n(K)$ the $n$th colored
Kauffman bracket, that is, the Kauffman bracket of the 
skein in the 3-dimensional sphere which is 
$K$ colored by the $n$th Jones-Wenzl idempotent. As such, $\kappa _1(K)$ is
the classical Kauffman bracket of the knot. There is a framing 
ambiguity in this definition, which we eliminate by always
working with the zero framing of the knot. 

\begin{theorem}
If a knot has the same noncommutative A-ideal
as the figure-eight knot, then all  colored Kauffman brackets are
the same.
\end{theorem}

\begin{proof} Two knots that have the same noncommutative A-ideals
have the same peripheral ideals (see \cite{frgelo}).
It has been shown in \cite{gelca2} that every element of 
the form $\sum_ic_i(p_i,q_i)_T$ in the peripheral ideal
gives rise to a family of relations of the form
\begin{eqnarray*}
& & \sum_ic_i(t^{(2n+p_i)q_i}[(-t^2)^{q_i}\kappa_{p_i+n}(K)
 -(-t^{-2})^{q_i}\kappa_{p_i+n-2}(K)]\\ & & +t^{-(2n-p_i)q_i
}[(-t^2)^{q_i}\kappa_{p_i-n}(K)-
(-t^{-2})^{q_i}\kappa_{p_i-n-2}(K)])=0.
\end{eqnarray*}
This is  the ``orthogonality relation between
the A-ideal and the Jones polynomial'' discovered in 
\cite{frgelo}. These relations indexed by the nonnegative
integer $n$ allow us to compute $\kappa _n(K)$ recursively. Of course,
one should note that $\kappa _0(K)=1$, $\kappa _{-1}(K)=0$, $\kappa _{-n-2}(K)=
-\kappa_n(K)$ for all $n$ and all knots.  

Let $K$ be the figure eight knot and
consider the element provided by Proposition \ref{peri31}, multiplied
by $-t^{-1}$ so that the coefficient of $(3,1)_T $ is $1$. 
The coefficient of $\kappa _3(K)$ in the  recursive
relation with $n=0$ 
is $-2t^{5}$, which is nonzero, so $\kappa _3(K)$ can be determined in
terms of $\kappa_2(K)$, $\kappa _1(K)$, $\kappa _0(K)=1$, $\kappa _{-1}(K)=0$
and $\kappa _{-2}(K)=-1$. For $n>1$, the coefficient of $\kappa _{n+3}(K)$
is $-t^{2n+5}$, which is always nonzero, so we can compute $\kappa _{n+3}(K)$,
$n\geq 0$, 
in terms of colored Kauffman brackets of lower color. 

On the other hand, an element of the form
\begin{eqnarray*}
\sum_q c_{2,q}(2,q)_T+\sum_qc_{1,q}(1,q)_T+\sum_qc_{0,q}(0,q)_T
\end{eqnarray*}
in the peripheral ideal gives rise to an orthogonality
relation of the form
\begin{eqnarray*}
\left[\sum_q(-1)^qc_{2,q}(t^{4q}+1)\right]\kappa _2(K)+
\left[\sum_q(-1)^qc_{1,q}(t^{3q}+t^q)\right]\kappa_1(K)+p(t)=0,
\end{eqnarray*}
where $p(t)$ is a Laurent polynomial in $t$.
 
Hence the elements provided by Proposition \ref{peri2}.
give rise to orthogonality relations of the form
\begin{eqnarray*}
& & P_1(t)\kappa_2(K)+Q_1(t)\kappa_1(K)=R_1(t)\\
& & P_2(t)\kappa _2(K)+Q_2(t)\kappa _1(K)=R_2(t).
\end{eqnarray*}
with $P_i(t),Q_i(t),R_i(t)$, $i=1,2$, Laurent polynomials in $t$. 
View this as a linear system of equations in the unknowns $\kappa_2(K)$
and $\kappa_1(K)$.

 One easily computes $P_1(t)=-t^{-6}+t^2$ and
$P_2(t)=-t^6+t^{-2}$, so $P_2(t)=(-t^4)P_1(t)$. 
Hence the system has zero determinant only if $Q_2(t)=(-t^4)Q_1(t)$. But
the term of highest degree in $Q_1(t)$ is $-t^{24}$ while  the 
term of highest degree in $Q_2(t)$ is $t^{16}$, so equality is
impossible. 

It follows that $\kappa_1(K)$ and $\kappa_2(K)$ are uniquely determined
by the two orthogonality relations, and the conclusion of the
theorem follows.

\end{proof}

\begin{cor}
Any knot having the same noncommutative A-ideal as the figure-eight knot
has the same Jones polynomial.
\end{cor}

In fact F. Nagasato wrote down the recursive relation obtained from the
first element of Proposition \ref{peri2} and checked that the recursion yields
Habiro's 
formulas for the colored  Kauffman brackets of the figure-eight
knot.

We have seen that the property stated in the above theorem
 is shared by the unknot \cite{gelca2} and
$(2,2p+1)$-torus knots \cite{gs}. There the property followed from 
 the existence of a polynomial of second degree in $l$ in the
noncommutative A-ideal. The present result  shows  that this
condition is not necessary. In \cite{gelca2} it  was also proved
that for any knot, the noncommutative A-ideal together with finitely many
colored Kauffman brackets determine all other colored Kauffman brackets
of the knot. All these facts suggest the following

\medskip

\begin{em}
{\bf Conjecture.} The noncommutative A-ideal of a knot determines all its
colored Kauffman brackets. In particular
if two knots have the same  noncommutative A-ideal
then they have the same Jones polynomial. 
\end{em}

\end{document}